\documentclass[b5paper,12pt] {article} 
\usepackage{buch} 
\begin{document}

    \title{Consistency Decision}

    \author{Michael Pfender\footnote{michael.pfender@alumni.tu-berlin.de}}
               
    \date{April 2014\\
         last revised \today}
    
\maketitle

\abstract{The consistency formula for set theory can be stated in terms
of the free-variables theory of primitive recursive maps. Free-variable
\pr\ predicates are decidable by set theory, main result here, built
on recursive evaluation of \pr\ map codes and soundness of that evaluation
in set theoretical frame: internal \pr\ map code equality is evaluated
into set theoretical equality. So the free-variable consistency 
predicate of set theory is decided by set theory, $\omega$-consistency
assumed. By Gödel's second incompleteness theorem on undecidability
of set theory's consistency formula by set theory under assumption
of this $\omega$-consistency, classical set theory turns out to be 
$\omega$-inconsistent.}

\tableofcontents

%\section*{Introduction} siehe abstract

\section{Primitive recursive maps}

Define the theory $\PR$ of objects and \pr\ maps as follows recursively 
as a subsystem of $\bf{set}$ theory $\T:$
\begin{itemize}
\item
the objects

$\one=\set{0},\N,\N\times\N,\ldots,A,\ldots,B,A\times B\ \etc$

\item
the map constants

$0:\one\to\N$ (\emph{zero}),\ $s=s(n)=n+1$ (\emph{successor}),
$\id_A:A \to A$ (\emph{identities}),
$\Pi:A\to \one$ (\emph{terminal} maps),
$l=l(a,b)=a:A\times B \to A,$
$r=r(a,b)=b:A\times B \to B$ (left and right \emph{projections});

\item
closure against (associative) map \emph{composition,}

$g\circ f = (g\circ f)(a) = g(f(a)):A\to B\to C;$

\item
closure against forming the \emph{induced} map 
$(f,g)=(f,g)(c)=(f(c),g(c)):C\to A\times B$ into a product, for given
components $f:C\to A,$ $g:C\to B,$ 

$l\circ (f,g)=f,\ r\circ (f,g)=g;$

\item
closure against forming the \emph{iterated} map
\begin{align*}
&f^\S=f^\S(a,n)=f^n(a):A\times\N\to A,\\
&f^0(a)=\id_A(a)=a,\\ 
&f^{sn}(a)=f^\S(a,sn)=(f\circ f^\S)(a,n)=f(f^n(a,n)).
\end{align*}  
\end{itemize}
Furthermore $\PR$ is to inherit from $\T$ uniqueness of the 
\emph{initialised iterated,} in order to inherit uniqueness
in the following \emph{full schema of primitive recursion:} 

\inference{(\mr{pr})}
{ $g=g(a):A\to B\ (\text{initialisation})$,\\
& $h=h((a,n),b):(A\times\N)\times B\to B\ (\text{step})$
}
{ $f=f(a,n):A\times\N\to B,$ \\
& $f(a,0)=g(a)$ \\
& $f(a,sn)=h((a,n),f(a))$ \\
& +\text{\emph{uniqueness} of such \pr\ defined map $f.$}
}
This schema allows in particular construction of \ul{for} loops,

\ul{for} $i:=1$ \ul{to} $n$ \ul{do}$\ldots$\ul{od}

as for verification if a given text (code) is an (arithmetised)
$\mathit{proof}$ of a given coded assertion, Gödel's \pr\
formula 45.\,$xBy,$ $x$ ist $\mathit{Beweis}$ von $y.$

(Formel 46.\,$\mathit{Bew\,y} = \exists x B y,$ $x$ is \emph{provable,} 
is not \pr)

\section{PR code sets and evaluation}

The \emph{map code set}---set of g\"odel numbers---we want 
to {\emph{evaluate}} is 
$\mrPR = \union_{A,B}[A,B] \subset \N$ in $\T,$ $ [A,B]=[A,B]_{\PR}$
the set of \pr\ map codes from $A$ to $B.$

Together with evaluation on suitable arguments it is recursively defined 
as follows:

\begin{itemize}
\item
Basic map constants $\mr{ba}$ in $\PR:$

\begin{itemize}
\item
$\code{0}\in[\one,\N]\subset\mrPR$ (zero),
$$ev(\code{0},0)=0,$$
$\code{s}\in[\N,\N]$ (successor), 
$$ev(\code{s},n)=s(n)=n+1,$$

\item
For an object $A$ $\code{\id_A}\in[A,A]$ (identity),
$$\ev(\code{\id_A},a)=\id_A(a)=a,$$
$\code{\Pi_A}\in[A,\one]$ (terminal map), 
$$\ev(\code{\Pi_A},a)=\Pi_A(a)=0.$$
%$\code{\Delta_A}\in[A,A\times A]$ (diagonal), 
%$$\ev(\code{\Delta_A},a)=\Delta_A(a)=(a,a).$$

\item
for objects $A,B$ $\code{l_{A,B}}\in[A\times B,A]$ (left projection),
$$\ev(\code{l_{A,B},(a,b)}=l_{A,B}(a,b)=a,$$ 
$\code{r_{A,B}}\in[A\times B,B]$ (right projection),
$$\ev(\code{r_{A,B},(a,b)}=r_{A,B}(a,b)=b.$$
\end{itemize}

\item
For $u\in[A,B],v\in[B,C]:$ $v\odot u \in [A,C]$ 

\quad(internal composition),
$$\ev(v\odot u,a)=\ev(v,\ev(u,a)).$$

\item
For $u\in[C,A],v\in[C,B]:$ $\an{u;v} \in [C,A\times B]$ 

\quad(induced map code into a product),
$$\ev(\an{u;v},c)=(\ev(u,c),\ev(v,c)).$$

\item
For $u\in [A,A]:$ $u^{\Dollar}\in[A\times\N,A]$ (iterated map code),
\begin{align*}
&\ev(u^{\Dollar},0)=\id_A(a)=a,\\
&\ev(u^{\Dollar},sn)=\ev(u,\ev(u^{\Dollar},n))\ \text{(double recursion)}
\end{align*}
This recursion \emph{terminates} in set theory $\T,$ with correct results:
\end{itemize}

\medskip
\textbf{Objectivity Theorem:} Evaluation $\ev$ is \emph{objective,} 
\ie\ for $f: A \to B$ in $\PR$ we have 
\begin{align*}
& \ev(\code{f},a) = f(a).
\end{align*}

\smallskip
\textbf{Proof} by substitution of codes of $\PR$ maps into code 
variables $u,v \in \mrPR \subset \N$ in the above double recursive
definition of evaluation, in particular:

\begin{itemize}
\item
composition
\begin{align*}
&\ev(\code{g} \odot \code{f},a) 
                        = \ev(\code{g},\ev(\code{f},a)), \\
&= g(f(a)) = (g \circ f)(a)
\end{align*}
recursively, and

\item
iteration
\begin{align*}
&\ev(\code{f}^{\Dollar},\an{a;sn}) 
      = \ev(\code{f},\ev(\code{f}^{\Dollar},\an{a;n}))\\
&= f(f^\S(a,n)) = f(f^n(a)) = f^{sn}(a)
\end{align*}
recursively.
\end{itemize}

\section{PR soundness within set theory}
 
Notion $f =^{\PR} g$ of \pr\ maps is externally \pr\ enumerated, 
by complexity of (binary) {deduction trees}. 

Internalising---\emph{formalising}---gives an internal notion 
of PR equality,
  $$u\,\check{=}_k\, v\in \mrPR \times \mrPR$$
coming by $k$th internal equation $\mathit{proved}$ by $k$th internal \emph{deduction tree} $\dtree_k.$ 

%which can be 
%canonically provided with arguments in $\X$---top down from (suitable) 
%argument $x$ given to the \emph{root} equation $u\,\check{=}_k\,v$ of 
%$\dtree_k.$

%We denote internal deduction tree argumented this way by 
%$\dtree_k/x,$ \emph{root} of $\dtree_k/x$ then is $u/x\,\check{=}_k\,v/x.$

\bigskip
\textbf{PR evaluation \emph{soundness} theorem framed by 
set theory $\T:$}
For \pr\ theory $\PR$ with its internal notion of equality `$\check{=}$' 
we have:
\begin{enumerate} [(i)]
\item
$\mrPR$ to $\T$ evaluation {soundness:}
\begin{align*}
\T \derives\ 
& u\,\check{=}\,v \implies \ev(u,x) = \ev(v,x) & (\bullet)
\end{align*}
Substituting in the above ``concrete'' $\PR$ codes into $u$ \resp\ $v,$ 
we get, by \emph{objectivity} of evaluation $\ev:$

\item $\T$-framed objective soundness of $\PR:$

For \pr maps $f,g: A \to B$
\begin{align*}
\T \derives\
& \code{f} \,\check{=}\, \code{g} \implies f(a) = g(a).
\end{align*}

\item 
Specialising to case $f :\,= \chi: A \to \two=\set{0,1}$ a \pr\ 
\emph{predicate,} and to $g :\,= \true,$ we get 
   
$\T$-framed \emph{logical soundness of $\PR:$}
\begin{align*}
& \T \derives\
      \exists k\,\Pro_{\PR} (k,\code{\chi}) 
                              \implies \forall\,x\,\chi(x):
\end{align*}
\emph{\textbf{If} a \pr\ predicate is---within $\T$---$\PR$-internally
\emph{provable,} then it holds in $\T$ for all of its arguments.}  
\end{enumerate}

\smallskip
\textbf{Proof} by
primitive recursion on $k,$ $\dtree_k$ the $k$\,th deduction tree of
the theory, \emph{proving} its root equation $u\,\check{=}_k\,v.$
These (argument-free) deduction trees are counted in 
lexicographical order. 

\smallskip
\textbf{Super Case} of \emph{equational} internal {axioms,} in particular
\begin{itemize}
\item
associativity of (internal) composition:

$\an{\an{w \odot v} \odot u}
            \,\check{=}\,\an{w \odot \an{v \odot u}} \implies$
\begin{align*}
& \ev\,(\an{w \odot v} \odot u,a) 
           = \ev\,(\an{w \odot v},\ev\,(u,a)) \\ 
& = \ev\,(w,\ev\,(v,\ev\,(u,a))) \\     
& = \ev\,(w,\ev\,(\an{v \odot u},a)) 
      = \ev\,(w \odot \an{v \odot u},a).
\end{align*}
This \textbf{proves} assertion $({\bullet})$ in present 
\emph{associativity-of-composition} case.

\item
Analogous \textbf{proof} for the other {flat,} equational cases,
namely \emph{reflexivity of equality,} \emph{left and right neutrality}
of identities, all substitution equations for the map 
constants, Godement's equations for the induced map:
$$l\odot\an{u;v}\,\check{=}\,u,\ r\odot\an{u;v}\,\check{=}\,v,$$ 
as well as \emph{surjective pairing} 
$$\an{l\odot w;r\odot w}\,\check{=}\,w$$
and distributivity equation
$$\an{u;v}\odot w\,\check{=}\,\an{u\odot w;v\odot w}$$
for composition with an induced. 

\item
\textbf{proof} of $(\bullet)$ 
for the last equational \textbf{case,} the    
 
\smallskip
\emph{Iteration step,}
    case of \emph{genuine iteration equation}
 
$u^{\Dollar} \odot \an{\code{\id} \# \code{s}
             \,\check{=}
                \,u \odot u^{\Dollar}},\ \#$ the internal cartesian product
                of map codes: 
\begin{align*}
\T \derives\
& \ev\,(u^{\Dollar} \odot 
    \an{\code{\id} \# \code{s}},\an{a;n}) & (1) \\
& = \ev\,(u^{\Dollar},\ev
           (\an{\code{\id} \# \code{s}},\an{a;n})) \\
%& = \ev\,(u^{\Dollar},\an{a;\code{s} \odot n}) \\
& = \ev\,(u^{\Dollar},\an{a;sn}) \\
%& = \ev\,(u^{[sn]},a) \\  
            %\quad \text{(by definition of $\ev$ step $e$)} 
%& = \ev\,(u \odot u^{[n]},a) \\
& = \ev\,(u,\ev(u^{\Dollar},\an{a;n}) \\
& = \ev\,(u \odot u^{\Dollar},\an{a;n}).  & (2)
\end{align*}

\end{itemize}

\textbf{Proof} of termination-conditioned inner soundness for the
remaining genuine \NAME{Horn} case axioms, of form 
$$u\,\check{=}_i\,v\,\land\,u'\,\check{=}_j\,v' 
  \implies w\,\check{=}_k\,w',\ i,j<k:$$

\smallskip
\textbf{Transitivity-of-equality} case
$$u\,\check{=}_i\,v\,\land\,v\,\check{=}_j\,w \implies u\,\check{=}_k\,w:$$
  
Evaluate at argument $a\in A$ and get in fact
\begin{align*} 
\T \derives\
& u\,\check{=}_k\,w \\ 
& \implies \ev(u,a) = \ev(v,a)\,\land\,\ev(v,a) = \ev(w,a) \\
& (\text{by hypothesis on}\ u,v) \\
& \implies \ev(u,a) = \ev(w,a): \\
& \text{transitivity export \qed\ in this case.}  
\end{align*}

%\medskip
%Case of \textbf{symmetry} axiom scheme for equality is obvious.

\medskip
\textbf{Compatibility case} of composition with 
equality,
\begin{align*} 
& u\,\checkeq\,u' \implies \an{v \odot u}\,\checkeq\,\an{v \odot u'}: \\
& \ev(v \odot u,a) = \ev(v,\ev(u,a)) = \ev(v,\ev(u',a)) \\
& = \ev(v \odot u',x),
\end{align*}
by hypothesis on $u\,\checkeq\,u'$ and by Leibniz' 
substitutivity in $\T,$ \qed\ in this first compatibility case.

%\smallskip
%Spread down arguments is more involved in  

\medskip
\textbf{Case} of composition with equality in second composition factor,
\begin{align*}
& v\,\checkeq_i\,v' \implies \an{v \odot u}\,\checkeq_k\,\an{v' \odot u}:\\
& \ev(\an{v \odot u},x) = \ev(v,\ev(u,x)) 
                          = \ev(v',\ev(u,x)) & (*) \\
& = \ev(\an{v'\odot u},x).
\end{align*} 
$(*)$ holds by $v\,\checkeq\,v',$ induction hypothesis on
$v,v',$ and Leibniz' substitutivity: same argument put into equal maps.

This proves soundness assertion $(\bullet)$ in this 2nd 
compatibility case.

\medskip
(Redundant) Case of \textbf{compatibility} of forming the induced map,
with equality, is analogous to compatibilities above, even
easier, since the two map codes concerned are 
independent from each other what concerns their domains.

\medskip
\textbf{(Final) Case} of Freyd's (internal) \textbf{uniqueness} 
of the \emph{initialised iterated,} is \textbf{case} 
\begin{align*}
& \an{w \odot \an{\code{\id};\code{0} \odot \code{\Pi}}\,\checkeq_i\,u} \\
& \land\,\an{w \odot \an{\code{\id} \# \code{s}}
                             \,\checkeq_j\,\an{v \odot w}} \\
& \implies w\,\checkeq_k\,v^{\Dollar} \odot \an{u \# \code{\id}}
\end{align*}
\textbf{Comment:} $w$ is here an internal \emph{comparison candidate} 
fullfilling the same internal \pr\ equations as the \emph{initialised iterated}
$\an{v^{\Dollar} \odot \an{u \# \code{\id}}}.$
It should be -- \textbf{is}: \emph{soundness} -- evaluated equal to the 
latter, on $A \times \N.$ 

\smallskip
Soundness \textbf{assertion} $(\bullet)$ for the 
present Freyd's \emph{uniqueness} \textbf{case} recurs
on $\checkeq_i,\ \checkeq_j$ turned into predicative equations
`$=$', these being already deduced, by hypothesis on $i,j < k.$ 
Further ingredients are transitivity of `$=$' and established 
properties of evaluation $\ev.$

\smallskip
So here is the remaining -- inductive -- \textbf{proof,} prepared by
\begin{align*}
\T \derives\ 
& \ev(w,\an{a;0}) = \ev(u;a) & (\bar{0}) \\  
& \qquad
    \text{as well as} \\   
& \ev(w,\an{a;sn}) = \ev\,(w \odot \an{\code{\id} \# \code{s}},\an{a;n}) \\ 
& = \ev\,(v \odot w,\an{a;n}), & (\bar{s})
\end{align*}
the same being true for 
  $w'  :\,= v^{\Dollar} \odot \an{u \# \code{\id}}$
in place of $w,$ once more by (characteristic) double recursive 
equations for $\ev,$ this time with respect to the 
\emph{initialised internal iterated} itself. 

\smallskip
$(\bar{0})$ and $(\bar{s})$ put together for both then show, by 
{induction} on \emph{iteration count} $n\in \N$---all other free 
variables $u,v,w,a$ together form the \emph{passive parameter} for this 
induction---\emph{soundness} assertion $({\bullet})$  for this 
\emph{Freyd's uniqueness} case, namely
\begin{align*}
\T \derives\ 
& \ev\,(w,\an{a;n}) =  
      \ev\,(v^{\Dollar} \odot \an{u \# \code{\id}},\an{a;n}).
\end{align*}
\textbf{Induction} runs as follows:

\textbf{Anchor} $n = 0:$

\smallskip 
$\ev\,(w,\an{a;0}) = \ev\,(u,a) = \ev\,(w' ,\an{a;0}),$
      
\medskip
\textbf{step:} 
\begin{align*}
& \ev(w,\an{a;n}) = \ev(w' ,\an{a;n}) \implies \\
& \ev\,(w,\an{a;sn}) = \ev(v,\ev(w,\an{a;n})) \\
& = \ev(v,\ev(w',\an{a;n})) 
    = \ev(w' ,\an{a;sn}), 
\end{align*}
%the latter since evaluation $\ev$ preserves predicative equality `$=$' (Leibniz)   
\textbf{\qed}\ 

%\medskip
%\textbf{Comment:}
%Already for stating the evaluations, we needed the---categorical, 
%free-variables theories $\PR, \PRa, \PRX, \PRXa$ of primitive recursion. %as %well as---for termination, even in classial frame $\T$---PR 
%%complexities within $\N[\omega].$ 
%Since this type of \textbf{soundness} 
%is a corner stone in our approach, the above complicated categorical 
%combinatorics seem to be necessary, even for the negative results on classical %foundations.

%Next section is the basis for the \emph{positive} logical Results
%concerning our \pr\ descent theory $\piR.$

\section{PR-predicate decision}
 
We consider here $\PR$ predicates for {decidability} by 
\textbf{set} theorie(s) $\T.$ Basic tool is 
$\T$-{framed soundness of} $\PR$ just above, namely
\inference{ } %({sound}\ \PR\backslash\T) }
{ $\chi = \chi(a): A \to \two$\ \,$\PR$ predicate }
{ $\T \derives\ \exists k\,\Pro_{\PR}(k,\code{\chi}) 
                                 \implies \forall a\,\chi(a).$ }

%This is an immediate consequence---$\T$ a (quantified) extension 
%of $\PRa$---of $\T$-dominated soundness of $\PRa$ just above.

\smallskip
Within $\T$ {define} for $\chi: A \to \two$ out of $\PR$ a 
partially defined (alleged, individual) $\mu$-recursive \emph{decision} 
$\nabla\chi: \one \parto \two$ by first fixing 
\emph{decision domain}
  $$D = D\chi :\,= \set{k \in \N:\neg\,\chi(\ct_A(k))
                           \,\lor\,\Pro_{\PR}(k,\code{\chi})},$$
$\ct_A: \N \to A$ (retractive) Cantor count of $A;$ and then,
with (partial) recursive $\mu D: \one \parto \N$ within $\T:$                           
\begin{align*}
& \nabla\chi \defeq
  \begin{cases}
  \false\ \myif\ \neg\,\chi(\ct_A(\mu D)) \\ 
    \quad (\emph{counterexample}), \\
  \true\ \myif\ \Pro_{\PR}(\mu D,\code{\chi}) \\
    \quad (\emph{internal proof}), \\
  \bot\ (\emph{undefined})\ \text{otherwise,}\ \ie\ \\
  \quad 
    \myif\ \forall a\,\chi(a)
      \,\land\,\forall k\,\neg\,\Pro_{\PR}(k,\code{\chi}).  
  \end{cases}
\end{align*}
[\,This (alleged) decision is apparently $\mu$-recursive within $\T,$
even if apriori only partially defined.]

\medskip
There is a first \emph{consistency} problem with this
\textbf{definition:} are the \emph{defined} cases \emph{disjoint?}

Yes, within frame theory $\T$ which \emph{soundly frames} 
theory $\PR:$ 
  $$\T \derives\ (\exists\,k \in \N)\,\Pro_{\PR}(k,\code{\chi}) 
                 \implies \forall a\,\chi(a).$$ 
We show now, that decision $\nabla\chi$ is \emph{totally defined,} the
undefined case does not arise, this for $\T$ $\omega$-consistent 
in Gödel's sense.

\medskip
We have the following complete -- metamathematical -- {case distinction} 
on $D=D_{\chi} \subseteq \N:$
\begin{itemize}
\item
\textbf{1st case,} termination: $D$ has at least one (``total'')
PR point $\one \to D \subseteq \N,$ and hence
  $$t = t_{\nabla\chi} \bydefeq \mu D = \min D: \one \to D$$ 
is a (total) \pr\ point. 

\textbf{Subcases:}
\begin{itemize}
\item
\textbf{1.1,} negative (total) \textbf{subcase:}
  
$\neg\chi\ct_A(t) = \true.$ 

[\,{Then} $\T \derives\ \nabla\chi = \false.$]
  
\item
\textbf{1.2,} positive (total) \textbf{subcase:}
  
$\ProPR(t,\code{\chi}) = \true.$ 

[\,{Then} $\T \derives\ \nabla\chi = \true,$ 

by $\T$-framed objective soundness of $\PR.$]
  
\smallskip
These two {subcases} are \emph{\textbf{disjoint,}} 
disjoint here by $\T$-framed soundness of theory $\PR$
which reads
\begin{align*}
\T \derives\ 
& \ProPR(k,\code{\chi}) \implies \forall a \chi(a),\ k\ \text{free,} 
\end{align*}
here in particular -- substitute $t: \one \to \N$ into $k$ free:
\begin{align*}
\piR \derives\ 
& \ProPR(t,\code{\chi}) \implies \forall a \chi(a). 
\end{align*}
 
So furthermore, by this framed soundness, in present \textbf{subcase:}
    $$\T \derives \forall a \chi(a)\,\land\,\ProPR(t,\code{\chi}).
                                                         \quad(\bullet)$$
\end{itemize}
  
\item
\textbf{2nd case,} derived non-termination: 

$T \derives\ D = \emptyset \identic \set{\N:\false_\N} \subset \N$

\qquad
[\,\text{then in particular}
     $\T \derives\ \forall a \neg\chi(a) = \false,$

\qquad
  \textbf{so} $\T \derives \forall a \chi(a)$ in this case\,],

{and furthermore}\ 
\begin{align*}
& \T \derives\ \forall k \neg\ProPR(k,\code{\chi}),\ \text{so}\\
&\T \derives\ \forall a \chi(a)
                \,\land\,\forall k \neg\ProPR(k,\code{\chi})\quad(*)
\end{align*}
in this case.

\item
\textbf{3rd,} remaining, \emph{ill} \textbf{case} is:

$D$ (metamathematically) \emph{has no (total) points $\one \to D,$ 
but is nevertheless not empty.} 
\end{itemize}

Take in the above the {(disjoint) union} of \textbf{2nd subcase} 
of \textbf{1st case,} $(\bullet),$ and of \textbf{2nd case,} 
$(*),$ as new case. And formalise last, remaining case. {Arrive at} the following 

\medskip
\textbf{Quasi-Decidability Theorem:} each \pr\ predicate $\chi: A \to \two$ 
gives rise within \textbf{set} theory $\T$ to the following
{complete (metamathematical) case distinction:}
\begin{enumerate}[(a)]
\item
$\T \derives \forall a \chi(a)$ {or else} 

\item
$\T \derives \neg\chi\ct_A t: \one \to D_{\chi} \to \two$

(\emph{defined counterexample}), {or else}

\item
$D = D_{\chi}$ \emph{non-empty, pointless,}
formally: in this \textbf{case} we would have within $\T:$ 
\begin{align*} 
& \T \derives\ \exists\,\hat{a} \in D, \\
& \text{{and} ``nevertheless'' {for each} \pr\ point}\ p: \one \to \N \\
& \T\derives\ p\not\in D. 
\end{align*} 
\end{enumerate}

We \textbf{rule out} the latter -- general -- possibility of a 
\emph{non-empty} predicate \emph{without \pr\ points}, for 
frame theory $\T$ by g\"odelian \emph{\textbf{assumption}} of
$\omega$-consistency. In fact it rules out above instance of
$\omega$-\emph{inconsistency:} all numerals $0,1,2,\ldots$
are \pr\ points. Hence it rules out -- in \emph{quasi-decidability} above -- possibility (c) for decision domain $D = D_{\chi} \subseteq \N$ 
of decision operator $\nabla$ for predicate $\chi: A \to \two,$ and we get

\medskip
\textbf{Decidability theorem:} Each free-variable \pr\ predicate 
$\chi: A \to \two$ gives rise to the following 
\textbf{complete case distinction} by \textbf{set} theory $\T:$

Under \emph{\textbf{assumption}} of $\omega$-consistency for $\T:$
\begin{itemize}
\item
$\T \derives\ \forall a \chi(a)$ (\emph{theorem}) {or} 

\item
$\T \derives (\exists\,a \in A)\,\neg\,\chi(a).$\ (\emph{counterexample})
\end{itemize}
Now take here for predicate $\chi,$ $\T$'s own free-variable \pr\ consistency 
formula
$$\ConT = \neg\,\ProT(k,\code{\false}): \N \to \two,$$ 
and get, under \emph{\textbf{assumption}} of $\omega$-consistency for $\T,$
a \textbf{consistency decision} $\nabla_{\Con_\T}$ for $\T$ by $\T.$

This contradiction to (the postcedent of) G\"odel's 
{2nd Incompleteness theorem} shows that the 
\emph{\textbf{assumption}} of $\omega$-consistency  
for \textbf{set} theories $\T$ must fail:

\medskip
\textbf{Set} theories $\T$ are $\omega$-inconsistent.

\smallskip
This concerns all classical \textbf{set} theories as in particular 
$\PM,\ \ZF,$ and $\NGB.$ 
%as well as already Peano arithmetic $\PA.$ 
%(NO??:Arithmetical Foundations) 
The reason is ubiquity of formal quantification within these (arithmetical) theories.

\textbf{Problem:} Does it concern Peano Arithmetic either?

\end{document}